\documentclass[final]{siamart0516}

\usepackage{amssymb}
\usepackage{mathabx,bm,xcolor,soul}
\newsiamthm{remark}{Remark}
\newsiamthm{example}{Example}

\numberwithin{equation}{section}

\newcommand{\R}{\mathbb{R}}

\newcommand{\id}{\ensuremath{\mathrm{id}}}
\newcommand{\tdd}{\ensuremath{\mathrm{tdd}}}
\newcommand{\<}{\langle}
\renewcommand{\>}{\rangle}

\DeclareMathOperator{\ddiv}{div}
\DeclareMathOperator{\sym}{sym}

\usepackage{mathtools}

\DeclarePairedDelimiterX{\inp}[2]{\langle}{\rangle}{#1, #2}

\begin{document} 
\title{A fractional-order trace-dev-div inequality}
\author{C.~Carstensen\thanks{Dep. of Maths,  Humboldt-Universit\"at zu Berlin, Germany; cc@math.hu-berlin.de}\and
N.~Heuer\thanks{Facultad de Matem\'aticas, Pontificia Universidad Cat\'olica de Chile,
Avenida Vicu\~na Mackenna 4860, Santiago,
Chile; nheuer@mat.uc.cl.
This author was supported by ANID-Chile through FONDECYT project 1230013.}}

\maketitle

\begin{small} {\bf Abstract.}
{The trace-dev-div inequality in $H^s$ controls the trace in the norm of $H^s$ by that of the deviatoric part plus the $H^{s-1}$ norm of 
the divergence of a quadratic tensor field different from the constant unit matrix. 
This is well known for $s=0$ and established for orders $0\le s\le 1$ and arbitrary space dimension in this note.
For mixed and least-squares finite element error analysis in linear elasticity, this inequality allows to establish robustness with respect to the Lam\'e parameter $\lambda$.}


\end{small}

\section{Main Result and Comments} 
This short note establishes the estimate 
\begin{equation}\label{cceqapx1}
 C_\tdd^{-1} \| \textrm{tr}\, \tau\|_{H^s(\Omega)}  
 \le \| \textrm{dev}\, \tau \|_{H^s(\Omega)} + \|  \textrm{div}\, \tau \|_{H^{s-1}(\Omega)} 
 \quad\text{for all } \tau\in \Sigma.
\end{equation}
This trace-deviator-divergence inequality (tr-dev-div for short) holds for any closed linear subspace $\Sigma$ of the Sobolev space $H^s(\Omega;\mathbb{R}^{n\times n})$
for $0\le s\le 1$ and a bounded Lipschitz domain $\Omega\subset \R^n$ with
arbitrary dimension $2\le n\in\mathbb{N}$
provided $\Sigma\not\ni \textrm{id}$. It is clear that the constant unit matrix $\textrm{id}$ violates  \eqref{cceqapx1} and has to be excluded;
but $\Sigma$ is closed and so an entire  neighbourhood 
(with respect to the norm in
$H^s(\Omega;\mathbb{R}^{n\times n})$)
of  $\textrm{id}$ lies outside  $\Sigma$. 

\begin{theorem}[tr-dev-div]\label{thmtr-dev-div}
If the closed linear subspace  $\Sigma$ of
$H^s(\Omega;\mathbb{R}^{n\times n})$ for  some $0\le s\le 1$
does not include  \id, then  \eqref{cceqapx1} holds for
a positive constant  $C_\tdd$. 
\end{theorem}

Some comments are in order before 
the proof of Theorem~\ref{thmtr-dev-div}  
follows below with  Bogovski\u\i's right-inverse of the 
divergence operator \cite{CostabelMacintosh,GHH}.

\begin{remark}[s=0]
The tr-dev-div inequality \eqref{cceqapx1}  is well established  in the least-squares 
and mixed finite element community for $s=0$ (when $H^0(\Omega)=L^2(\Omega)$)
\cite{BBF,CaiS}.
It has been utilised for Stokes equations or linear elasticity 
to provide $\lambda$-robust error control.
\end{remark}

\begin{remark}[symmetric variant]
Linear elasticity involves symmetric tensors 
$\mathbb{S}=\sym\mathbb{R}^{n\times n}$ 
and $\Sigma$ may well be a 
closed linear subspace of $H^s(\Omega;\mathbb{S})$
with $\id\notin\Sigma$
(because $H^s(\Omega;\mathbb{S})$ is a closed subspace of 
$H^s(\Omega;\mathbb{R}^{n\times n})$).
\end{remark}

\begin{remark}[s=1, n=2]
The tr-dev-div inequality   \eqref{cceqapx1}   leads to  the $\lambda$-robust
regularity estimate \cite[Eq (11.2.33)]{BS-3.ed} from \cite[Thm A.1]{Vogelius} of 
$\lambda\|  \ddiv u\|_{H^s(\Omega)}$ for 
pure homogeneous Dirichlet or Neumann boundary conditions 
on a convex polygon, $n=2$, $s=1$, and a constant Lam\'e parameter $\lambda>0$.
\end{remark}

\begin{remark}[$\lambda$-robust regularity]\label{remarklambdarobustregularity}
The paper \cite{Roessle_00_CSR} computes the singular functions for rather general (mixed) boundary
conditions in $n=2$ dimensions and provides $\varepsilon(u)$ in $H^s(\Omega)$ for 
explicit $s>0$. Given positive material constants $\lambda$ and $\mu$, called the
Lam\'e constants, the elastic stress reads 
$  \sigma=2\mu\varepsilon{u}+\lambda(\ddiv u)\id$ and satisfies 
$\ddiv\sigma=-f\in L^ 2(\Omega)^2$ for the given source term in 
$L^ 2(\Omega)^2$. Therefore, 
 the tr-dev-div inequality   \eqref{cceqapx1} 
provides $\lambda\,  \ddiv u\in H^s(\Omega)$ and 
controls its norm $\lambda\|\ddiv u\|_{H^s(\Omega)}$ 
with $\lambda$-independent constants and the $H^{s-1}$ norm of the source term. 
\end{remark}

\begin{example}[convex domain and pure boundary condition for $n=2,3$]
Given a source term $f\in L^ 2(\Omega)^n$ in a  convex bounded domain 
$\Omega$ and suppose either pure Dirichlet conditions 
(resp.  pure Neumann conditions) are prescribed with data 
$u_D \in H^{3/2}(\partial\Omega)^ n$ (resp. compatible data 
 $g\in H^{1/2}(\partial\Omega)^ n$
with $\int_\Omega f\, dx+\int_{\partial\Omega} g\, ds =0$).
Then the weak solution $u\in H^1(\Omega)^n$ and the associated 
stress $\sigma =2\mu\varepsilon{u}+\lambda(\ddiv u)\id\in L^ 2(\Omega;\mathbb{S}) $ 
satisfy  $u\in H^2(\Omega)^n$ and $\sigma \in H^ 1(\Omega;\mathbb{S})$ 
and  
\[
\| u\|_{H^ 2(\Omega)}+\lambda\|\ddiv u\|_{H^1(\Omega)}\le C_\mathrm{reg} 
( \|f\|_{L^ 2(\Omega)}+\| u_D\|_{H^{3/2}(\Omega)})
\]
for the pure Dirichlet problem 
(resp. $\le C_\mathrm{reg}  ( \|f\|_{L^ 2(\Omega)}+\| g\|_{H^{1/2}(\Omega)})$ for the pure
compatible Neumann problem).
The point is that the constant $C_\mathrm{reg} $ exclusively depends on $\Omega$ and $\mu$
but {\em not}  on $\lambda$.
(The proof follows from well known higher regularity for convex domains and 
Remark~\ref{remarklambdarobustregularity}.)
\end{example}

\section{Sobolev spaces by interpolation} 
Let $\widetilde H^m(\Omega):= H^m_0(\Omega):= \overline{\mathcal{D}(\Omega)} $ be the  closed subspace of  $H^m(\Omega) $ 
for all non-negative integers  $m\in\mathbb{N}_0$. The Sobolev spaces for the 
negative integers follow by duality (with pivot space $L^2$ and pairing $\<\cdot,\cdot\>$) 
\[
H^{-m}(\Omega):= \widetilde H^m(\Omega)^*\quad\text{and}\quad \widetilde  H^{-m}(\Omega):= H^m(\Omega)^*\quad\text{for } m \in\mathbb{N}_0.
\]
Given the  Sobolev spaces $\widetilde H^m(\Omega)$ and  $H^m(\Omega)$ for all integers $m\in\mathbb{Z}$, define 
\[
H^r(\Omega)=[ H^{m}(\Omega), H^{m+1}(\Omega)]_s
\quad\text{and}\quad
\widetilde H^r(\Omega)=[\widetilde H^{m}(\Omega) ,\widetilde H^{m+1}(\Omega)]_s
\quad\text{for } r=m+s
\] 
with $0<s<1$ and $m\in\mathbb{Z}$ by complex interpolation 
\cite{Bergh-Loefstroem,Calderon_63_ISI,Lions,Tartar,Triebel}.
A real interpolation method could also be utilised and may lead to additional equivalence constants. The complex interpolation,  after
Calder\'on \cite{Calderon_63_ISI} and Lions \cite{Lions},
allows for equal norms on the dual of an interpolation 
space compared to the interpolation of the duals. 
We frequently utilise the symmetry $[X,Y]_s=[Y,X]_{1-s}$ and 
the duality theorem
 $([X,Y]_s)^*=[X^*,Y^*]_{s}$
with equal norms for reflexive spaces $X,Y$,
cf.~\cite[Theorems~4.2.1(a),~4.3.1]{Bergh-Loefstroem} 
and \cite[Corollary~4.5.2]{Bergh-Loefstroem}, respectively.

 We underline that $H^s(\Omega)=\widetilde H^s(\Omega)$ for 
 $0\le s<1/2$ while, for $1/2\le s\le 1$, $\widetilde H^s(\Omega)$ 
 is a strict  subspace of  $H^s(\Omega) $ (with vanishing traces for $1/2<s\le 1$). 
 The definition of  $\widetilde H^{1/2}(\Omega)$ by interpolation
results in the Lions--Magenes space $H^{1/2}_{00}(\Omega)$, a strict subset of $H_0^{1/2}(\Omega)=\overline{\mathcal{D}(\Omega)}$
\cite{Grisvard}.
For instance,   $H^{-1/2}(\Omega)$  is  dual to  $\widetilde H^{1/2}(\Omega)\equiv H^{1/2}_{00}(\Omega)$. 
The norm in  $H^{-1/2}(\Omega)$ arises in \eqref{cceqapx1} for $s=1/2$ as the dual norm to 
 $\widetilde H^{1/2}(\Omega)\equiv H^{1/2}_{00}(\Omega)$
which has a stronger than $ H^{1/2}(\Omega)$. Hence 
\eqref{cceqapx1} holds in particular when  $\|  \textrm{div}\, \tau \|_{H^{-1/2}(\Omega)} $ is confused with
  $\|  \textrm{div}\, \tau \|_{\widetilde H^{-1/2}(\Omega)}
  = \|  \textrm{div}\, \tau \|_{H^{1/2}_0(\Omega)^*}\ge  \|  \textrm{div}\, \tau \|_{H^{-1/2}(\Omega)}$.

\section{Gradient and divergence in $\widetilde H^r(\Omega)$}
 The gradient $\nabla : H^m(\Omega)\to H^{m-1}(\Omega)^n$ and its restriction  
 $\widetilde \nabla := \nabla |_{\widetilde H^m(\Omega)}  : \widetilde H^{m}(\Omega)\to  \widetilde H^{m-1}(\Omega)^n$ as well as 
 the divergence  $\textrm{div} : H^m(\Omega)^n\to H^{m-1}(\Omega)$   and its restriction   
 $\widetilde{\textrm{div}}:= \textrm{div} |_{\widetilde H^m(\Omega)^n}  : \widetilde H^{m}(\Omega)^n \to  \widetilde H^{m-1}(\Omega)$
are well defined pointwise a.e. for any positive integer $m\in\mathbb{N}$. 
For $m\in\mathbb{N}$,  those definitions provide linear and bounded operators 
 \begin{align*}
 \widetilde \nabla \in & \hspace{1mm} L(\widetilde  H^m(\Omega);  \widetilde  H^{m-1}(\Omega)^n)&\text{and} & &
 \widetilde{\textrm{div}}\in & \hspace{1mm} L(\widetilde H^{m}(\Omega)^n ;  \widetilde H^{m-1}(\Omega)), \\
  \nabla \in & \hspace{1mm} L(H^m(\Omega);  H^{m-1}(\Omega)^n)  &\text{with dual} & &
   \nabla^* \in & \hspace{1mm} L(\widetilde H^{1-m}(\Omega)^n ;\widetilde H^{-m}(\Omega) ),
  \\
\textrm{div} \in & \hspace{1mm} L(H^m(\Omega)^n; H^{m-1}(\Omega)) &\text{with dual} & &\textrm{div}^*  \in 
& \hspace{1mm} L(\widetilde H^{1-m}(\Omega) ; \widetilde H^{-m}(\Omega)^n) .
\end{align*}
Integration by parts for 
$(\varphi,\Psi)\in 
\mathcal{D}(\Omega)\times C^\infty_{\textrm{c}}(\overline{\Omega}) ^n$,
\[
\int_\Omega \varphi\, \textrm{div}\, \Psi\, dx+\int_\Omega  \Psi\cdot\nabla \varphi\, dx =0 ,
\]
and the  density of $ C^\infty_{\textrm{c}}(\overline{\Omega})$ in $H^r(\Omega)$ 
and $\mathcal{D}(\Omega)$ in $\widetilde H^r(\Omega)$ justify the definition of 
\[
\widetilde \nabla:= -   \textrm{div}^* \in L( \widetilde H^{1-m}(\Omega); \widetilde H^{-m}(\Omega)^n) 
\]
for $m\in\mathbb{N}$.
Integration by parts for $(\varphi,\Psi)\in  C^\infty_{\textrm{c}}(\overline{\Omega}) \times \mathcal{D}(\Omega)^n$ justifies the
definition of  $\widetilde{\textrm{div}}:=- \nabla ^* \in  L(\widetilde H^{1-m}(\Omega)^n ;\widetilde H^{-m}(\Omega) )$ for $m\in\mathbb{N}$. 
In summary this defines 
\begin{equation}\label{cceqapx2}
\widetilde \nabla \in L( \widetilde H^{r}(\Omega); \widetilde H^{r-1}(\Omega)^n) 
\quad\text{and}\quad 
\widetilde{\textrm{div}}\in  L(\widetilde H^{r}(\Omega)^n ; \widetilde H^{r-1}(\Omega) )
\end{equation}
for any integer  $r\in\mathbb{Z}$ and thereafter complex interpolation defines the operators \eqref{cceqapx2} 
for any real $r$. 
The norms of the resulting operators are controlled by interpolation: The operator norm of 
$\widetilde \nabla$ and $\widetilde{\textrm{div}}$ from  \eqref{cceqapx2}  is at most 
$\sqrt{n}$ and $1$, respectively,  for $r=0$; while  $1$ and $\sqrt{n}$, respectively,  for  $r=1$.
For  the operator norms and $0\le s\le 1$,  
\begin{equation}\label{cceqapx3}
\| \widetilde \nabla \|_{ L( \widetilde H^{1-s}(\Omega); \widetilde H^{-s}(\Omega)^n) }\le n^{s/2} 
\quad\text{and}\quad 
\| \widetilde{\textrm{div}}\|_{  L(\widetilde H^{1-s}(\Omega)^n ; \widetilde H^{-s}(\Omega) )}\le  n^{(1-s)/2} 
\end{equation}
follows by interpolation.

\begin{remark}[alternative Sobolev spaces]
The extension of  
 $\nabla : H^s(\Omega)\to H^{s-1}(\Omega)^n$ for a general real  $s$ is linear and bounded except for $s=1/2$ 
in the standard textbook  \cite[Thm 1.4.4.6]{Grisvard}. This is not a contradiction to \eqref{cceqapx3} for all real $s$
because \cite{Grisvard} 
utilises
$H^{-1/2}(\Omega)$ as the topological dual of
$H^{1/2}_0(\Omega)=H^{1/2}(\Omega)$ which gives $\widetilde H^{-1/2}(\Omega)$.
\end{remark}

\section{Poincar\'e inequality in $H^s(\Omega)\cap L^2_0(\Omega)$ for $0\le s\le 1$}
This subsection investigates the constant $\|B\|_s$ in the estimate 
\begin{equation}\label{cceqapx4} 
\|B\|_s^{-1}\, \|g\|_{H^s(\Omega)}\le \sup_{v\in \mathcal{D}(\Omega)^n\setminus\{0\}} 
\frac{ \int_\Omega g\, \textrm{div} \, v\, dx}{ \| v\|_{ \widetilde H^{1-s}(\Omega)} }
\le  n^{(1-s)/2} \, \|g\|_{H^s(\Omega)}
\end{equation}
for all  $g\in H^s(\Omega)\cap L^2_0(\Omega)$. Here and throughout this paper $L^2_0(\Omega)=L^2(\Omega)/\R$
abbreviates the Lebesgue functions in $L^2(\Omega)$ with integral zero over the domain.
Since the supremum defines the norm $\| \nabla g\|_{ H^{s-1}(\Omega)} $,
\eqref{cceqapx4}  is a Poincar\'e inequality (sometimes called Poincar\'e--Friedrichs inequality)  for $s=0$ and an important duality estimate 
for $s=1$ usually attributed to Ne\v cas with applications to a proof of Korn's  inequality. 

The remaining parts of this subsection are devoted to the proof of 
\eqref{cceqapx4} based on
Bogovski\u\i's right-inverse of the divergence operator. 
The proof departs from duality
\begin{equation}\label{cceqapx5}  
\|g\|_{H^s(\Omega)}=\sup_{f\in L^2_0(\Omega)\setminus\{0\}} 
\frac{ \int_\Omega g\, f\, dx}{ \| f \|_{ \widetilde H^{-s}(\Omega)} }
\quad\textrm{for all } g\in H^s(\Omega)\cap L^2_0(\Omega).
\end{equation}   
Since $1\in H^{s}(\Omega)$, the integral mean $\overline{f}=\<f,1\>/|\Omega|$ is well defined for $f\in  \widetilde H^{-s}(\Omega)$  
and the continuous linear operator
$f\mapsto f-\overline{f}$ has a norm at most one in  $L(L^2(\Omega))$ and in  $L(\widetilde H^{-1}(\Omega))$ and so also in 
$L(\widetilde H^{-s}(\Omega))$, i.e., 
$  \| f-\overline{f} \|_{ \widetilde H^{-s}(\Omega)} \le  \| f \|_{ \widetilde H^{-s}(\Omega)} $ for $f\in L^2(\Omega)$.
This provides the inequality in the dual norm 
\[
 \|g\|_{H^s(\Omega)}= \sup_{ f\in  \widetilde H^{-s}(\Omega) \setminus\{0\}} \frac{ \<f,g\> }{  \| f \|_{ \widetilde H^{-s}(\Omega)} }
 \le  \sup_{ f\in  \widetilde H^{-s}(\Omega) \setminus\{0\}} \frac{ \<f-\overline{f},g\> }{  \| f - \overline{f} \|_{ \widetilde H^{-s}(\Omega)} }
\]
with $\<\overline{f},g\>=0$ in the last step. The upper bound allows the substitution of $ f - \overline{f}$ for $f\in  \widetilde H^{-s}(\Omega)$ by 
$f\in  \widetilde H^{-s}(\Omega)\cap L^2_0(\Omega) $ and therefore provides ``$\le$'' in  \eqref{cceqapx5}; its converse ``$\ge$'' is obvious.

Bogovski\u\i's right-inverse of the divergence operator \cite{CostabelMacintosh} serves as  a standard tool in the analysis of the Stokes equations 
and \cite[Thm 2.5]{GHH} provides the existence of $B\in L( \widetilde H^{-s}(\Omega); \widetilde H^{1-s}(\Omega)^n$ with operator norm $\|B\|_s$
and the right-inverse property 
\[
v:= B f  \in H^1_0(\Omega) ^n \quad\textrm{satisfies } \textrm{div} \, v = f\quad\textrm{for all }f\in L^2_0(\Omega).
\]
This enables  the substitution of $f \in L^2_0(\Omega)$ in the supremum of \eqref{cceqapx5} by $ \textrm{div} \, v = f$ for 
$v:= B f  \in H^1_0(\Omega) ^n$. This and  $\| v \|_{ \widetilde H^{1-s}(\Omega)}\le \| B\|_s\,  \| f \|_{ \widetilde H^{-s}(\Omega)}$ reveals the 
first inequality in  \eqref{cceqapx4}. The second follows from 
$ \int_\Omega g\, \textrm{div} \, v\, dx\le  \| \widetilde{\textrm{div}} \, v\|_{\widetilde H^{-s}(\Omega)} \|g\|_{H^s(\Omega)} $ and 
$ \|\widetilde{\textrm{div}} \, v\|_{\widetilde H^{-s}(\Omega)}\le    n^{(1-s)/2} \|  v\|_{\widetilde H^{1-s}(\Omega)}$ by \eqref{cceqapx3}.

\section{Proof of  tr-dev-div  \eqref{cceqapx1} for
a particular $\Sigma$}
We select $\tau\in H^s_\mathrm{tr0}(\Omega;\mathbb{R}^{n\times n}):=$
$\{\sigma\in H^{s}(\Omega;\mathbb{R}^{n\times n}):
\; \mathrm{tr}\,\sigma\in L^2_0(\Omega)\}$,
a closed subspace of $H^{s}(\Omega;\mathbb{R}^{n\times n})$. Then
$g:=\textrm{tr} \,\tau \in L^2_0(\Omega)$ and 
\eqref{cceqapx4} motivates the analysis of $ \int_\Omega g\, \textrm{div} \, v\, dx$ 
for some $v\in \mathcal{D}(\Omega)^n$. Since $n^{-1} \textrm{tr} \,\tau \, \id= \tau -\textrm{dev} \, \tau $ a.e. in $\Omega$, some algebra 
reveals 
\[
n^{-1} \int_\Omega  \textrm{tr} \,\tau  \, \textrm{div} \, v\, dx = \int_\Omega  \tau  : D  v\, dx -  \int_\Omega  (\textrm{dev}\, \tau ) : D  v\, dx
\]
with the functional matrix $Dv\in L^2(\Omega;\R^{n\times n}) $ that consists of the row-wise action of the gradient $\widetilde \nabla$.
 Recall
 $v\in \mathcal{D}(\Omega)^n$ and infer 
 \[
 \int_\Omega  \tau  : D  v\, dx=-\<\textrm{div}\,\tau,   v\>\; \le \; \| \textrm{div}\,\tau\|_{ H^{s-1}(\Omega)}   \| v\|_{ \widetilde H^{1-s}(\Omega)} .
\]
Since $D v$ is the row-wise application of $\widetilde \nabla$, \eqref{cceqapx3}
 provides $\| Dv \|_{ \widetilde H^{-s}(\Omega)} \le n^{s/2}\,  \| v\|_{ \widetilde H^{1-s}(\Omega)} $.  Thus
\[ 
  -  \int_\Omega  (\textrm{dev}\, \tau ) : D  v\, dx \le \|\textrm{dev}\, \tau \|_{ H^{s}(\Omega) } \| Dv \|_{ \widetilde H^{-s}(\Omega)}
  \le n^{s/2}\,  \|\textrm{dev}\, \tau \|_{ H^{s}(\Omega) } \| v\|_{ \widetilde H^{1-s}(\Omega)}.
 \]
 Altogether, for any   $v\in \mathcal{D}(\Omega)^n$ with $\| v\|_{ \widetilde H^{1-s}(\Omega)}=1$, 
\[
n^{-1} \int_\Omega  \textrm{tr} \,\tau  \, \textrm{div} \, v\, dx 
\le n^{s/2}\,  \|\textrm{dev}\, \tau \|_{ H^{s}(\Omega) }
+  \| \textrm{div}\,\tau\|_{ H^{s-1}(\Omega)} .
\] 
This and the first estimate in \eqref{cceqapx4} conclude the proof of 
 \eqref{cceqapx1} for
$\Sigma=H^{s}_\mathrm{tr0}(\Omega;\mathbb{R}^{n\times n})$
with the constant  $ C_\tdd=n^{1+s/2} \|B\|_s $.
 
\section{Proof of tr-dev-div \eqref{cceqapx1} for general \mbox{$\Sigma$}}
Recall $\Sigma$ is a closed linear subspace of
$H^{s}(\Omega;\mathbb{R}^{n\times n})$
that does not include the identity $\id$. 
Suppose the   trace-dev-div estimate  \eqref{cceqapx1}  fails and so we find some sequence
$\tau_1,\tau_2,\tau_3,\dots $ in
$\Sigma\subset H^{s}(\Omega;\mathbb{R}^{n\times n})$
with (for all $j\in\mathbb{N}$)  
\[
 \| \textrm{tr} \,\tau_j\|_{H^s(\Omega)}=1\quad\text{and}\quad  \lim_{j\to\infty}  \| \textrm{dev} \,\tau_j\|_{H^s(\Omega)}=
 0= \lim_{j\to\infty}  \| \textrm{div} \,\tau_j\|_{H^{s-1}(\Omega)}.
\]
Without loss of generality suppose that the integral mean $\mu_j:=|\Omega|^{-1}\int_\Omega \textrm{tr} \,\tau_j\, dx\ge 0$ 
is non-negative. Owing to the previous subsection,  
$\dot\tau_j:=\tau_j - \mu_j n^{-1} \id \in$
$H^{s}_\mathrm{tr0}(\Omega;\mathbb{R}^{n\times n})$
satisfies  \eqref{cceqapx1}.
Since the divergence and the deviator part  of $\tau_j$ and $\dot\tau_j$ are identical, we infer
 \[
 \lim_{j\to\infty}\|  \textrm{tr} \,\dot \tau_j\|_{H^s(\Omega)} =0 < 1 =  \lim_{j\to\infty}\|  \textrm{tr} \, \tau_j\|_{H^s(\Omega)},
 \]
 whence $ \lim_{j\to\infty}  \mu_j \| 1\|_{H^s(\Omega)}= 1$. In other words, 
 $ \lim_{j\to\infty} \tau_j =  \| 1\|_{H^s(\Omega)}^{-1} n^{-1}\id$. Since 
 $\Sigma\ni\tau_j $ is closed, it 
contains the limit 
 $ \| 1\|_{H^s(\Omega)}^{-1} n^{-1}\id\in \Sigma$ 
 and the vector space $\Sigma$ contains $\id$.

\bibliographystyle{amsplain}
\bibliography{abstractpaper.bib}
\end{document}